\documentclass[a4paper,12pt]{article}
\usepackage{amsmath}
\usepackage{amsthm}
\usepackage{amssymb}
\usepackage{amscd}
 \makeatletter

 \@addtoreset{equation}{section}
  \makeatother

\theoremstyle{plain}
\newtheorem{thm}{Theorem}[section]
\newtheorem{prop}[thm]{Proposition}

\newtheorem{lem}[thm]{Lemma}
\theoremstyle{definition}
\newtheorem{defn}{Definition}[section]

\theoremstyle{remark}

\begin{document}
\title{Extensions of Current Groups on \(S^3\) and the Adjoint Representations}
\author{ Tosiaki Kori\\
Department of Mathematics\\
School of Science and Engineering\\
 Waseda University \\
3-4-1 Okubo, Shinjuku-ku
Tokyo, Japan.\\ e-mail: kori@waseda.jp
}

\date{ }
\maketitle

\begin{abstract} 
Let \(\Omega^3(SU(n))\) be  the Lie group of based mappings from \(S^3\) to  \(SU(n)\).   We construct a Lie group extension of \(\Omega^3(SU(n))\) for \(n\geq 3\) by the abelian group  \(\exp 2\pi i\,{\cal A}_3^{\ast}\,\),  where  \({\cal A}_3^{\ast}\) is the affine dual of   the space  of \(SU(n)\)-connections  on \(S^3\).     
     J. Mickelsson in 1987 constructed a similar Lie group extension.    In this article we give several improvement of  his results, especially we give a precise description of the extension of those components that are not the identity component,.   We also correct several argument about the extension of \(\Omega^3(SU(2))\)  which seems not to be exact in Mickelsson's work, though his observation about the fact that  the extension of \(\Omega^3(SU(2))\)  reduces to the extension by \({\bf Z}_2\) is correct.     Then we shall  investigate the adjoint representation of the  Lie group extension of  \(\Omega^3(SU(n))\) for \(n\geq 3\).
\end{abstract}

\section{Introduction} 

Let \(G\) be a compact Lie group and let \(MG=Map(M, G)\) be the set of smooth mappings from a manifold \(M\) to  \(G\) that are based at some point of \(M\).   Groups of the form \(MG\) have  been a subject of investigation both from a purely mathematical standpoint and from quantum field theory.   In quantum field theory they appear as current groups or gauge transformation groups.    In the simplest case \(M\) is the unit circle \(S^1\) and \(LG=Map(S^1,G)\) is a loop group.   Loop groups and their representation theory have been fully worked out.   \(LG\) turned out to behave like a compact Lie group and the highly developed theory of finite dimensional Lie groups was extended to the infinite dimensional group \(LG\).    \(LG\) appears in the simplified model of quantum field theory where the space is one-dimensional and many important facts in the representation theory of loop groups were first discovered by physicists.   It turned out that in many applications to field theory one must deal with certain extensions of loop groups and their associated Lie algebras.    The central extension of \(Map(S^1,\,Lie\,G)\) is an affine Kac-Moody algebra and the highest weight theory of finite dimensional Lie algebra was extended to this case.   \cite{B},  \cite{Kh-W},  \cite{PS} and \cite{Se} are good references to study these subjects.   But we know little about the generalization of the above picture to higher dimensional space \(M\).       
 In 1987 J. Mickelsson~\cite{Mic} gave a Lie group extension of 
\(\Omega^3G=Map(S^3, G)\) for  \(G=SU(n)\).     
Recently the author in \cite{K2} constructed the pre-quantization of the moduli space of flat connections on a four-manifold.   The group \(\Omega^3G\) acts symplectically on this moduli space, but it does not lift to an action on the pre-quantization bundle.   He showed that Mickelsson's extension is necessary to lift this action.   So Mickelsson's extension is an appropriate and natural extension of  \(\Omega^3G\).    In the following we shall explain it for the case \(n\geq 3\).     Let \(f\in S^4G\) and let \(\mathbf {f}\in D^5G\) be the extension of \(f\) to the 5-dimensional disk.    Such an extension is possible since \(\pi_4(SU(n))\) vanishes for \(n\geq  3\).    
We consider the five dimensional mapping degree ( or the five dimensional Chern-Simons form restricted to the pure gauges ):
\begin{equation}C_5(f)=\frac{i}{240\pi^3}\int_{D^5}\,tr( d\mathbf {f}\cdot \mathbf{f}^{-1})^5.
\end{equation}
Since \(H^5(SU(n),\mathbf{ Z})=\mathbf{ Z}\),  \(C_5(f)\) is defined by modulus \(\mathbf{ Z}\) independently of the extension \(\mathbf {f}\).    
It holds that 
\begin{equation}
C_5(fg)-C_5(f)-C_5(g)=\beta_{S^4}(f,g)\quad\mod \mathbf{ Z} .\label{ibeta}
\end{equation}
Where
\begin{eqnarray}
 \beta_{S^4}(f,g)&=&\frac{i}{48\pi^3}\int_{S^4}\,c^{2,1}(f,g).  \label{ibeta2}
 \end{eqnarray}
for a 4-form valued 2-cocycle \(c^{2,1}\) on the Lie group \(S^4G\).
Then we see that 
\begin{eqnarray}
\chi_{S^4} (f,g)&=& \exp 2\pi i\left(\,C_5(g)-C_5(f) \,\right) \\[0.2cm]
&=& \exp 2\pi i\left(\,\beta_{S^4}(f,f^{-1}g)+C_5(f^{-1}g) \,\right) 
\end{eqnarray}
satisfies the cocycle condition ; 
\(
\chi_{S^4} (f,g)\chi_{S^4} (g,h)=\chi_{S^4} (f,h).
\)
 Let   \(\Omega^3_0G\) be the connected component of the identity in \(\Omega^3G\).      Now we consider the group \({\cal G}_0\) of gauge transformations on the hemisphere \(D^4\) that are identity on the boundary \(S^3\);   \({\cal G}_0=\{g\in D^4G;\,g\vert S^3=1\}\).   Then \(\Omega^3_0G\simeq D^4G/{\cal G}_0\).     
 We observe that the definition of \(\chi_{S^4}(f,g)\) is extended to those mappings \(f,g\in D^4G\) such that \(f\vert S^3=g\vert S^3\).    In fact extend \(f^{-1}g\in{\cal G}_0\) by \(1\) on \(S^4\setminus D^4\), then the integral on the right-hand side of ( \ref{ibeta2}) is done over \(D^4\) .     We define  
 \begin{equation}
 \chi_{D^4}(f,g)=\exp 2\pi i\left(\,\beta_{D^4}(f,f^{-1}g)+C_5(f^{-1}g\vee 1') \,\right) .
\end{equation}
\(\chi_{D^4} (f,g)
\)
satisfies the cocycle condition.   
 Then we have a line bundle \(\,L=D^4G\times \mathbf{C}/{\cal G}_0\longrightarrow \Omega^3_0G\,\).          Contrary to the case of loop groups,  \(L\setminus\{0\}\) has no group structure.   
 Instead Mickelsson~\cite{Mic} considered the associated principal bundle;  
\begin{equation}
\pi:\,\widehat{\Omega_0 G}=L\times_{{\cal G}_0}Map({\cal A}_3,U(1))\longrightarrow \Omega^3_0G.
\end{equation}  
\({\cal A}_3\) being the space of connections on \(S^3\), and gave the group structure on it by the 2-cocycle known as   
{\it Mickelsson's 2-cocycle} (\ref{Mic3}).     
Thus he got a group extension  of \(\Omega^3_0G\) by the abelian group \(Map({\cal A}_3,U(1))\).     To have the group extension of the total space \(\Omega^3G\) we consider, instead of \(D^4G\), a mapping cone \(TG\) of the set of smooth mappings from \(T=S^3\times [0,1]\) to \(G\).       In~\cite{Mi, Mic} the abelian extension of \(\Omega^3_0G\) is fully evolved but as for \(\Omega^3G\) was only sketched the outline.     
So we shall modify in section 2 the Mickelsson's argument to fit to \(TG\).     There is another distinguished  improvement of Mickelsson's theory.   Instead of \(Map({\cal A}_3,U(1))\) we consider a more tight  subspace \(\exp 2\pi i{\cal A}_3^{\ast}\,\) of  \(Map({\cal A}_3,U(1))\), where \({\cal A}^{\ast}_3\) is the affine dual of \({\cal A}_3\), that is, the vector space of the maps \(l:{\cal A}_3\longrightarrow \mathbf{R}\)  that satisfies \(l(pA_1+(1-p)A_2)=pl(A_1)+(1-p)l(A_2)\)  for all  \(A_1,\,A_2\in{\cal A}_3\) and \(p\in \mathbf{ R}\).    \(\Omega^3G\) acts on \({\cal A}^{\ast}_3\) and    Mickelsson's 2-cocycle belongs to \({\cal A}^{\ast}_3\).     So we obtain the extension; 
\begin{equation}
1 \longrightarrow \exp 2\pi i{\cal A}_3^{\ast}\longrightarrow \widehat{\Omega G}\longrightarrow \Omega^3G\longrightarrow 1.
\end{equation}
Here we mention the group extension of  \(\Omega^3(SU(2))\).    In this case both \(\beta_{D^4}\) and 
\(C_5\) vanish and we have only the trivial extension of \(\Omega^3_0(SU(2))\).        
The argument in~\cite{Mic} to yield this fact contains several misunderstandings and we shall give a correct  proof.     As for the component of the group \(\Omega^3(SU(2))\) other than the connected component of the identity,  we shall find an interesting phenomenon that  concerns Witten's fermionization principle \cite{W} and this is due to the fact \(\pi_4(SU(2))={\bf Z}_2\).   Witten proved the formula:
  \begin{equation}
 \exp 2\pi iC_5(\widetilde g)=\epsilon(g),\qquad \widetilde g=\left(\begin{array}{cc}g&0\\0&1\end{array}\right),\quad  g\in S^4(SU(2)), \label{iC5}
 \end{equation}
 where \(\epsilon(g)=\pm 1\) if \(g\) represents the trivial ( respectively non-trivial ) homotopy class of \(\pi_4(SU(2))\).    
 If we embed \(SU(2)\) in \(SU(3)\) and if we consider the 
 restriction of the group extension \(\widehat{\Omega (SU(3))}\) of \(\Omega^3(SU(3))\) to the embedded subgroup  \(\Omega^3(SU(2))\) we have an extension of \(\Omega^3(SU(2))\) by  \(\mathbf{ Z}_2\) with the transition function given by \(\chi(f,g)=\exp 2\pi iC_5(\widetilde f^{-1}\,\widetilde g)=\epsilon(f^{-1}g)\).     This extension is not topologically trivial but turns out to be algebraically trivial.
 
     In section 3 we shall discuss the corresponding Lie algebra extensions and give the formula of  adjoint representation of \(\widehat {\Omega  (SU(n))}\), for \(n\geq 3\).

\section{  Basic properties on current groups}

\subsection{Descent equations}

Let  \(G=SU(n)\).   Let \(N\) be an oriented 5- manifold.    As typical examples we are thinking of  the 5-sphere \(N=S^5\),  the 5-dimensional disk \(N=D^5\) and  \(N=S^3\times D^2\) where \(D^2\) is the 2-dimensional disk.   Let 
 \(P=N\times G\) be the trivial \(G\)-principal bundle over \(N\).   Let \({\cal A}(N)\) denote the space of connections on \(N\).     The group of gauge transformations on \(N\) is denoted by \({\cal G}(N)\).   Since \(P\) is a trivial principal bundle \({\cal G}(N)\) is the space \(Map(N,G)\) of smooth mappings from \(N\) to \(G\) that are pointed at some point.

  Let \(\Omega^q(N)\) be the differential \(q\)-forms on \(N\) and let \(V^q\) be the vector space of polynomials \(\Phi=\Phi(A)\) of \(A\in {\cal A}(N)\)  and its curvature \(F_A\) that take values in \(\Omega^q(N)\).  The curvature \(F_A\) of a connection \(A\) will be often abbreviated to \(F\).   The group of gauge transformations \({\cal G}\) acts on \(V^q\) by \((g\cdot \Phi)(A)=\Phi(g^{-1}\cdot A)\).   We shall investigate the double complex 
\[C^{p,q}=C^p({\cal G}, V^{q+3}),\] 
that is doubly graded by the chain degree \(p\) and the differential form
degree \(q\).   Let \(\,d:C^{p,q}\longrightarrow C^{p,q+1}\) be the exterior differentiation.    The coboundary operator \(\delta: C^{p,q}\longrightarrow C^{p+1,q}\) is given by
\begin{equation*}
\begin{split}
(\delta\, c^p)(g_1,g_2,\cdots,g_{p+1})&=
g_1\cdot c^p(g_2,\cdots,g_{p+1})+(-1)^{p+1}c^p(g_1,g_2,\cdots,g_p)\\[0.2cm]
&+\sum_{k=1}^p(-1)^kc^p(g_1,\cdots,g_{k-1},g_kg_{k+1},g_{k+2},\cdots,g_{p+1}).
\end{split}
\end{equation*}
We introduce the following cochains:
\begin{eqnarray*}
c^{0,2}(A)&=&\,tr\,(AF^2-\frac12
A^3F+\frac1{10}A^5),\\[0.2cm]
c^{1,2}(g)&=& c^{0,2}(dg\,g^{-1})\,=\,\frac1{10}
tr(dg\,g^{-1})^5,\\[0.2cm]
c^{1,1}(g; A)&=& tr[\,-\frac12V(AF+FA-A^3)+
\frac14(VA)^2+\frac12V^3A\,],\\[0.2cm]
&&\qquad\mbox{where \(V=dg\,g^{-1}\)}, \\[0.2cm]
c^{2,1}(g_1,g_2)&=&c^{1,1}(g_2 ;\,g_1^{-1}dg_1\,),\\[0.2cm]
c^{2,0}(g_1,g_2;A)&=&\frac12 tr\lbrack\, \left (g_1^{-1}dg_1\,dg_2g_2^{-1}
- dg_2g_2^{-1}\,
g_1^{-1}dg_1\right )g_1^{-1}Ag_1\,\rbrack ,\\[0.2cm]
c^{3,0}(g_1,g_2,g_3)&=&c^{2,0}(g_2,g_3,g_1^{-1}dg_1).
\end{eqnarray*}
In the above  \(dg\,g^{-1}\) is the 1-form on \(N\) that is the pullback by \(g\in Map(N,G)\) of the Maurer -Cartan form.    The curvature \(F_A\) is abbreviated to \(F\).

\begin{prop}\cite{K2}\label{WZdescent}~
The cochains \(c^{p,q}\in C^{p,q}\), \(0\leq p,q\leq 3\),  satisfy the relations:
\begin{eqnarray}
dc^{p,3-p}+(-1)^p\delta c^{p-1,3-p+1}&=&0\label{WZdescent1}\\[0.2cm]
dc^{p,2-p}+(-1)^{p}\delta c^{p-1,3-p}&=&-c^{p,3-p}\label{WZdescent2}\\[0.2cm]
c^{0,3}=0,\qquad c^{p,q}&=&0\qquad \mbox{if }\,p+q\neq 2,3 \,.\nonumber
\end{eqnarray}
\end{prop}

The {\it Chern-Simons form}  on \(N\) is by definition.   
\begin{equation}
c^{0,2}(A)=tr(\,A F^2-\frac12 A^3F+\frac1{10} A^5\,),\qquad  A\in{\cal A}(N)\quad F=F_A.
\end{equation}  
\begin{prop}
The variation of the Chern-Simons form along the \({\cal G}(N)\)-orbit is given by:  
\begin{equation}
c^{0,2}( g\cdot  A)-c^{0,2}( A)=d\, c^{1,1}( g,\, A)+c^{1,2}( g), \quad g\in {\cal G}(N),
\label{eq:variation} 
\end{equation}
\end{prop}
This follows from (\ref{WZdescent2}).

\vspace{0.3cm}

In the following we shall show that when we consider the Lie group \(SU(2)\)  the above quantities \(c^{1,2}\), \(c^{2,1}\), \(c^{2.0}\) and \(c^{3,0}\) vanish.    

\begin{lem}
Let \(\alpha,\,\beta,\,\gamma\) be 1-forms on a 3-manifold valued in the Lie algebra \(su(2)\).   Then 
\begin{equation}
tr[\,(\alpha\,\beta-\beta\,\alpha)\,\gamma\,]=0.
\end{equation}
\end{lem}

{\it Proof}

Let \(dx^i;\,i=1,2,3\) be the local coordinates and let \(e_a;\,a=1,2,3\) be the basis of \(su(2)\);
\[e_ae_b=-e_be_a=-\epsilon_{abc}e_c,\qquad (e_a)^2=-I,\]
where \(\epsilon\) is totally antisymmetric in \(a,b,c\) and \(\epsilon_{123}=1\).   Let 
\[\alpha=\sum_i\,\alpha_idx^i=\sum_i(\sum_a\alpha^a_ie_a)dx^i,\]
and similarly for \(\beta=\sum_j\,\beta_jdx^j\) and \(\gamma=\sum_k\,\gamma_kdx^k\).
Then we have 
\[(\alpha\beta-\beta\alpha)\gamma=\left(\sum\,\epsilon_{ijk}(\alpha_i\beta_j+\beta_j\alpha_i)\gamma_k\right)dx^1dx^2dx^3.\]
Since \[\alpha_i\beta_j+\beta_j\alpha_i=\sum_{a,b}\,\alpha_i^a\beta_j^b(e_ae_b+e_be_a)=(2\sum_a\alpha^a_i\beta^a_j)I,\]
we have 
\[(\alpha\beta-\beta\alpha)\gamma=\sum\,\epsilon_{ijk}C_{ij}\otimes \gamma_k\,,\]
where \(C_{ij}=2\sum_a\alpha^a_i\beta^a_j\,dx^1dx^2dx^3\).    The trace of the last 3-form is \(0\).   \hfill\qed

\begin{prop}\label{nul}
For \(G=SU(2)\), we have 
\begin{equation}
c^{2,0}=c^{3,0}=0\qquad c^{2,1}=0 . \label{su2}
\end{equation}
\end{prop}
 
 From the previous lemma it follows that  \(c^{2,0}=c^{3,0}=0\) on any three-manifold.     Now  \(c^{2,1}\) is given by 
\[c^{2,1}(g_1,g_2)=tr[\,\frac12VA^3+
\frac14(VA)^2+\frac12V^3A\,],\]
with \(A=g_1^{-1}dg_1\) and \(V=dg_2\,g_2^{-1}\).    For any \(su(2)\) valued 1-form \(\alpha=\sum_a\alpha^ae_a\) we have 
\[\alpha^3=-\left( \sum\,\epsilon_{a b c}\alpha^a\alpha^b\alpha^c\right)I\,,\]
where \(e_a;\,a=1,2,3\) are the basis of \(su(2)\) and \(\alpha^a\) are 1-forms.     Then, for any \(su(2)\) valued 1-form \(\beta\), 
\[\alpha^3\beta=\left( \sum\,\epsilon_{a b c}\alpha^a\alpha^b\alpha^c\right)\sum_p \beta^pe_p\,.\]
Hence  \(tr[\alpha^3\beta]=0\).   This yields the vanishing of \(tr[VA^3]\) and 
\(tr[V^3A\,]\).    
Therefore, if we let 
\(V=\sum V_idx^i=\sum V^a_ie_adx^i\) and \(A=\sum A_jdx^j=\sum A^b_je_bdx^j\),  we have 
\begin{eqnarray*}
c^{2,1}(g_1,g_2)&=&\frac14 tr[\,(VA)^2]= \frac14\sum \epsilon^{ijkl}tr[\,V_iA_jV_kA_l]\\[0.2cm]
&=&\frac12\sum\epsilon^{ijkl}(V^a_iA^a_j)(V^p_kA^p_l)-\frac14\sum\,\epsilon^{ijkl}tr[\,\epsilon^{abc}V^a_iA^b_je_c\,\epsilon^{pqr}V^p_kA^q_le_r\,]\\[0.2cm]
&=&\frac12\sum_{a\neq p}\epsilon^{ijkl}(V^a_iA^a_j)(V^p_kA^p_l)\,+\, \frac12\sum_{a\neq b}\epsilon^{ijkl}\,\,V^a_iA^b_jV^b_kA^a_l\\[0.2cm]
&=&\frac12\sum_{a\neq b}\epsilon^{ijkl}(V^a_iA^a_j)(V^b_kA^b_l)\,+\, \frac12\sum_{a\neq b}\epsilon^{ijkl}\,(V^a_iA^a_l)(V^b_kA^b_j)=0\,.\end{eqnarray*}
\hfill\qed

\subsection{Descent equations for current algebras}

In this paragraph we suppose that \(G=SU(n)\), \(n\geq 3\).   
We shall study the descent equations for the Lie algebra of infinitesimal gauge transformations \(Lie\,
{\cal G}\) .    We consider the double complex
\[
E^{p,q}=C^p(Lie\,{\cal G},V^{q+3})\,,
\]
that is doubly graded by the chain degree \(p\) and the differential form degree \(q\).   The infinitesimal action of \(\xi\in Lie\,{\cal G}\) on \(V^q\) is given by 
\((\xi\cdot\Phi)(A)=\frac{d}{dt}\vert_{t=0}\Phi(\exp(-t\xi)\cdot A)=\Phi(-d_A\xi)\).
 The coboundary operator \(\delta: E^{p,q}\longrightarrow E^{p+1,q}\) is defined by
\begin{equation*}
\begin{split}
(\delta\, e^p)(\xi_1,\xi_2,\cdots,\xi_{p+1})&=
\sum_{i<j}\,(-1)^{i+j}e^p([\xi_i,\xi_j],\xi_1,\cdots,\hat\xi_i,\cdots,\hat\xi_j,\cdots,\xi_{p+1})  \\[0.2cm]
&+\sum_{k=1}^{p+1}(-1)^{k+1}\xi_k\cdot e^p(\xi_1,\cdots,\xi_{k-1},\xi_{k+1},\xi_{k+2},\cdots,\xi_{p+1}).
\end{split}
\end{equation*}
We put
\begin{eqnarray*}
e^{1,1}(\xi;A)&=&\frac{d}{dt}\lvert_{t=0}\, c^{1,1}(\exp t\xi;\,A),\\[0.2cm]
 e^{2,0}(\xi,\eta,A)&=&\, \frac{d}{ds}\lvert_{s=0}\frac{d}{dt}\lvert_{t=0}\,
\,c^{2,0}(\exp s\xi,\exp t\eta\,;\,A).
 \end{eqnarray*}
Then we have
\begin{prop}
\begin{eqnarray}
\delta e^{1,1}&=&-de^{2,0} ,\\[0.2cm]
e^{1,1}(\xi;A)&=& tr[\,\frac12(AF+FA-A^3)d\xi\,],\nonumber \\[0.2cm]
 e^{2,0}(\xi,\eta\,;A)&=&\frac12\, tr[\,(\,d\xi d\eta-d\eta d\xi\,)A\,]\nonumber
\end{eqnarray}
\end{prop}
In fact from the definition we have 
\[(\delta e^{1,1})(\xi,\eta;A)=\xi\cdot e^{1,1}(\eta;A)-\eta\cdot e^{1,1}(\xi;A)-e^{1,1}([\xi,\eta];A).\]
Then the formulas
\begin{eqnarray*}
\xi\cdot A=d_A\xi&=&d\xi+A\xi-\xi A,\\[0.2cm]
 \xi\cdot (dA)&=&[dA,\xi]-[A,d\xi]=dA\xi-\xi dA-Ad\xi-d\xi A
 \end{eqnarray*}
 yield the desired equations. \hfill\qed

Let \(S^3(Lie\,G)\) be the Lie algebra of the based mappings rom \(S^3\) to \(Lie\,G\). 

 We put
\begin{eqnarray}
\omega(\,\xi,\eta;A\,)\,&=&\, - \frac{1}{12\pi^3}\int_{S^3}\,e^{2,0}(\xi,\eta;A)\nonumber\\[0.2cm]
&=&\, - \frac{1}{24\pi^3}\int_{S^3}\,tr[\,(d\xi d \eta-d\eta d\xi)A\,],
\end{eqnarray}
for \(
  A\in{\cal A}_3\) and \(\xi,\eta\in S^3Lie\,G\), 
and denote 
 \begin{equation}\label{omega}
 \omega_f(\xi,\eta)\,=\,\omega(\xi,\eta\,;\,f^{-1}df\,).
\end{equation}
\begin{prop}
\(\omega_f\)   is a closed 2-form on \(\Omega^3G\), hence \((\Omega^3G,\omega)\) is a pre-symplectic space.
\end{prop}
In fact, 
the exterior differential \(\widetilde d\,\omega_f\) of \(\omega_f\) on \(\Omega^3_0G\) becomes 
\begin{eqnarray*}
(\widetilde  d\,\omega_f)(\xi,\eta,\zeta)&=&\,
\frac{d}{dt}\vert_{t=0}\,\omega_{\exp t\zeta}(\xi,\eta)
\\[0.2cm]
&=&
-\frac{1}{24\pi^3}\int_{S^3}\,d\,tr[(d\xi d \eta-d\eta d\xi)\zeta]=0.
\end{eqnarray*}
 \hfill\qed

\subsection{Basic properties}

\subsubsection{  }

Let \(M\) be a compact four-manifold possibly with non-empty boundary \(\partial M\).    Let \(G=SU(n)\), \(n\geq 3\).   
 In the following we write by \(MG\) the set of smooth mappings $f$ from \(M\) to \(G\) that
are based at some point \(p_0\in M\).       

\begin{defn}
For \(f,\,g\in MG\) we put
\begin{equation}
 \beta_{M}(f,g)=\frac{i}{24\pi^3}\int_{M}\,c^{2,1}(f,g) \label{Mic2}
 \end{equation}
 \end{defn}
 \begin{defn}[Mickelsson's 2-cocycle]\cite{Mi, Mic}   
 For \(f,\,g\in MG\) we put
 \begin{eqnarray}
 \gamma_{M}(f,g\,; A)&=&
\frac{i}{24\pi^3}\int_{M}(\delta c^{1,1})(\,f,\,g\,;\,A) \label{Mic3}
\\[0.2cm]
&=&- \frac{i}{24\pi^3}\int_{\partial{M}}c^{2,0}(f,g\,;\,A)+\frac{i}{24\pi^3}\int_{M}\,c^{2,1}(f,g)\,,\nonumber\\[0.3cm]
&=&- \frac{i}{24\pi^3}\int_{\partial{M}}c^{2,0}(f,g\,;\,A)+\beta_{M}(f,g) 
\end{eqnarray}
\end{defn}
 If \(\partial M=\emptyset\),   we have
\begin{equation}
\beta_{M}(f,g)=\gamma_{M}(f,g;A).
\end{equation}

\begin{lem}
Let  \(f,g,h\in MG\).    
  We have 
\begin{equation}
\gamma_{M}(f,g;A)+\gamma_{M}(fg,h;A)=\gamma_{M}(g,h;A)+\gamma_{M}(f,gh;A)\,.\label{gamma}
\end{equation}
If moreover either \(\partial M=\emptyset\) or at least one of 
\(f,g,h\) is constant on \(\partial M\),   then 
\begin{equation}
\beta_{M}(f,g)+\beta_{M}(fg,h)=\beta_{M}(g,h)+\beta_{M}(f,gh)\,.\label{beta}
\end{equation}
\end{lem}
{\it Proof}

The definition (\ref{Mic3}) of \(\gamma_M\) implies that  \(\gamma_M\) is a coboundary, hence it satisfies the cocycle property (\ref{gamma}).     Next the relation 
\( \delta c^{2,1}=dc^{3,0}\), (\ref{WZdescent1}), implies 
 \[
 \beta_M(g,h)- \beta_M(f,g)- \beta_M(fg,h)+ \beta_M(f,gh)=\frac{i}{24\pi^3}\int_{\partial{M}}
 c^{3,0}(f,g,h).\]
  From the formula of \(c^{3,0}\) we see that the right hand side vanishes if  \(\partial M=\emptyset\) or at least one of 
\(f,g,h\) is constant on \(\partial M\), that implies (\ref {beta}).     \hfill\qed

\begin{lem}
\begin{eqnarray}
 \beta_{M}(f,f^{-1})&=&\gamma_{M}(f,f^{-1};A)=0,\\[0,2cm]
\beta_{M}(fg,g^{-1})&=&-\beta_{M}(f,\,g)\,=\,\beta_M(f^{-1},\,fg).\label{beta2}
\end{eqnarray}
\end{lem}
In fact we have 
\begin{eqnarray}
c^{2,1}(f,f^{-1})&=&c^{2,0}(f,f^{-1};A)=0,\nonumber\\[0.2cm]
c^{2,1}(fg,g^{-1})&=&-\,c^{2,1}(f,g)\,=\, c^{2,1}(f^{-1},fg).\nonumber
\end{eqnarray}
All these relations follows from Proposition~\ref{WZdescent} by direct calculation.

\subsubsection{Polyakov-Wiegmann formula for \(SU(n)\) with \(n\geq 3\) }

Now we suppose that \(G=SU(n)\) with \(n\geq 3\) and that \(M\) is a 4 dimensional manifold that is the boundary of a 5-dimensional connected simply connected manifold \(N\); \(\partial N=M\).   For example \(M=S^4=\partial D^5\) or \(M=S^3\times S^1=\partial ( S^3\times D^2)\).    Since \(\pi_4G=1\)  every \(g\in MG\) has an extension \({\bf g}\in NG\).    For \(g\in MG\)  we define 
\begin{equation}\label{C5}
C_5(g)=\frac{i}{24\pi^3}\int_{N}\, c^{1,2}({\bf g})\,= \frac{i}{240\pi^3}\int_{N} tr(d{\bf g}\cdot
{{\bf g}}^{-1})^5 .
\end{equation}
 \(C_5(g)\) may depend on the extension 
 but it can be shown that the difference 
of two extensions is an integer because of \(H^5(G, {\bf Z})={\bf Z}\).   Hence  \(C_5(g)\) is well defined\(\mod{\bf Z}\), or 
\(\exp(2\pi i C_5(g)\,)$ is well defined independently of the
extension.
\begin{lem}[Polyakov-Wiegmann]\label{Poly}~\cite{Mi, P}~  For \(f,\,g\in
M(SU(n))\), \(n\geq 3\), we have
\begin{equation}
C_5(fg)=C_5(f)+C_5(g)+ \beta_M(f,g)\quad
\mod{\bf Z}.\end{equation}
\end{lem}
From  (\ref{WZdescent1}),   
\(\delta c^{1,2}=-dc^{2.1}\).      
Integration over \(N\) proves the lemma.
 \hfill\qed
 
  \begin{lem}\label{MiC5}~\cite{Mic}~~Let  \(T=S^3\times [0,1]\).   Let \(f,\,g \in
T(SU(n))\),  \(n\geq 3\).    Suppose that  \(g(\cdot,0)=g(\cdot,1)=1\).   Then we have
\begin{equation}
C_5(fgf^{-1})=C_5(g)+\beta_T(fg, f^{-1})+\beta_T(f,g)\quad \mod{\bf Z}.
\end{equation}
\end{lem}

\subsection{\(2\pi\) rotation of \(\Omega^3(SU(2))\) in \(\Omega^3(SU(3))\)}
 
  Here we shall study the case for \(SU(2)\).    Since \(c^{2,1}=0\) and \(c^{2,0}=0\) from  Proposition \ref{nul},  \(\beta_M\) and \(\gamma_M\) in the preceding subsection can play no role.    But \(c^{1,2}\), hence \(C_5(g)\), will be an important quantity.   Since \(\pi_4(SU(2))={\bf Z}_2\),  \(g\in M(SU(2))\) does not necessarily have an extension to its five-dimensional counterpart as in the previous subsection, so we have no definition of  \(C_5(g)\).   But when we embed  \(M(SU(2))\) in  \(M(SU(3))\) this quantity may be defined and it represents actually a quantity 
that reflect the fact \(\pi_4(SU(2))={\bf Z}_2\).   

Let \(M=S^3\times S^1\).   We look on  \(M\) as the boundary of the five dimensional manifold \(Q=S^3\times D^2\),  where  \(D^2\) is the two-dimensional disk.   
 \(M(SU(2))\) is considered as a subgroup of  \(M(SU(3))\) by the embedding 
\begin{equation}
M(SU(2))\ni u\longrightarrow \widetilde u=\left(\begin{array}{rl}u&0\\0&1\end{array}\right) \in M(SU(3)).
\label{embed}
\end{equation}
      Then the functional 
 \(C_5(\widetilde u)\) is well defined modulo \({\bf Z}\):
 \[C_5(\widetilde u)=\frac{i}{240\pi^3}\int_{Q} tr(d {\bf u}\cdot
{ {\bf u}}^{-1})^5 ,\]
where \(\widetilde u\in M(SU(3))\) is extended to \({\bf u}\in Q(SU(3))\).   
  In \cite{W} Witten showed that  \(C_5(\widetilde u)\)  depends only on the homotopy class of \(u\in \pi_4(SU(2))\).     \(C_5(\widetilde u)=0\mod {\bf Z}\) if \( u\) is in the trivial homotopy class in \(\pi_4(SU(2))\).      On the other hand  \(C_5(\widetilde u)=-\frac12\mod {\bf Z}\) for \(u\)  in the non-trivial homotopy class in \(\pi_4(SU(2))\),  \cite{W}.    
  
 For \(f\in \Omega^3(SU(2))\),  Witten investigated the process of a \(2\pi\) rotation of 
\(\widetilde f=\left(\begin{array}{rl}f&0\\[0.1cm] 0&1\end{array}\right)  \) 
inside \(\Omega^3(SU(3))\).
   The path \(\{\,{\widetilde u}_f(t)\,\}_{0\leq t\leq 1}\) in \(\Omega^3(SU(3))\) obtained by rotating  \(\widetilde f\)  by a \(2\pi\) angle is chosen to be 
\begin{eqnarray}
{\widetilde u}_f(x,t)&=&\left(\begin{array}{rcl}e^{\pi it}&0&0\\[0.2cm]0&e^{-\pi it}&0\\[0.2cm]0&0&1
\end{array}\right)\,\widetilde f(x)\,\left(\begin{array}{rcl}e^{-\pi it}&0&0\\[0.2cm]0&e^{ \pi it}&0\\[0.2cm]0&0&1
\end{array}\right)\nonumber\\[0.2cm]
&=&\left(\begin{array}{rcl}1&0&0\\[0.2cm]0&e^{-2\pi it}&0\\[0.2cm]0&0&e^{ 2\pi it}
\end{array}\right)\,\widetilde f(x)\,\left(\begin{array}{rcl}1&0&0\\[0.2cm]0&e^{ 2\pi it}&0\\[0.2cm]0&0&e^{- 2\pi it}
\end{array}\right).\label{rotate}
\end{eqnarray}
Then  \({\widetilde u}_f\in M(SU(3))\).     We extend  \({\widetilde u}_f\) to \(Q=S^3\times D^2\)  by 
\[{\bf u}_f(x,t,r)=a(r,t)\,\widetilde f(x)\, a(r,t)^{-1},\]
where
\[ a(r,t)=\left(\begin{array}{rcl}1&0&0\\[0.2cm]0&re^{-2\pi it}&\sqrt{1-r^2}\\[0.2cm]0&-\sqrt{1-r^2}&re^{2\pi it}
\end{array}\right), \quad  0\leq r\leq 1.\]
By the first form of (\ref{rotate}) we see that if \(f\in \Omega^3_0(SU(2))\), that is, if \(f\)  is homotopic to the identity of \(SU(2)\), then \(C_5({\widetilde u}_f)=0\).     
While if \(\widetilde f\) is an instanton, that is, if \(f\in \Omega^3(SU(2))\) and \(\deg f=1\), then the second form of  (\ref{rotate}) describes the \(2\pi\) rotation of instanton and we find 
\[C_5({\widetilde u}_f)=\frac{i}{240\pi^3}\int_{Q} tr(d{\bf u}_f\cdot
{\bf u}_f^{-1}\,)^5=-\frac12\qquad\mod{\bf Z}.\] 
Thus for the non-contractible path \({\widetilde u}_f\) in \(\Omega^3(SU(2))\) we have \(C_5({\widetilde u}_f)    
=-\frac12\).
Since \(\pi_4((SU(2))=\pi_1(\Omega^3(SU(2))\) the non-trivial homotopy class in \(\pi_4((SU(2))\) corresponds to the homotopy class of non-contractible paths in \(\Omega^3(SU(2))\).   Therefore we have \(C_5(\widetilde u)=-\frac12\)\(\mod {\bf Z}\) for 
 \(u\)  in the non-trivial homotopy class in \(\pi_4(SU(2))\) .   
\begin{defn}\label{epsilon}
Let  \(M= S^3\times S^1\).
We put, for \(u\in M(SU(2))\), 
 \begin{equation*}
  \epsilon(u)=\left\{\begin{array}{cc}  1 ,&\qquad\mbox{if \(u\) is in the trivial homotopy class of \(\pi_4(SU(2))\),} \\[0.2cm]
  -1 ,&\qquad\mbox{if \(u\) is in the non-trivial homotopy class of \(\pi_4(SU(2))\).}
  \end{array}\right.
  \end{equation*}
  \end{defn}
\begin{prop}
\label{Cepsil}
Let \(u\in M(SU(2))\).   Then 
\(\exp 2\pi i C_5(\widetilde u)\) depends only on the homotopy class of \(u\) in \(\pi_4(SU(2))\) and 
\begin{equation}
\epsilon(u)=\exp 2\pi i C_5(\widetilde u).
\end{equation}
\end{prop}  
\begin{lem}
\label{PWU2}
Let \(u,\,v\in M(SU(2))\).      We have the following  formulas.
\begin{eqnarray}
C_5(\,\widetilde u\,\widetilde v\,)&=&C_5(\,\widetilde u\, )+C_5(\,\widetilde v\, ) \quad \mod{\bf Z},\\[0.3cm]
 \epsilon (uv)&=&\epsilon(u)\epsilon(v) ,\\[0.2cm]
 \epsilon(uvu^{-1})&=&\epsilon(v).
\end{eqnarray}
\end{lem}
{\it Proof}

 The product \(uv\) of  \(u,\,v \in M(SU(2))\)
is in the trivial homotopy class if both are in the same homotopy class of \(\pi_4(SU(2))\), and  \(uv\) is in the non-trivial homotopy class if \([u]\) and \([v]\) are in the distinct homotopy class of \(\pi_4(SU(2))\).    
So we have the relation  \( \epsilon (uv)=\epsilon(u)\epsilon(v)\).   From Proposition \ref{Cepsil} we have 
\(C_5(\,\widetilde u\,\widetilde v\,)=C_5(\,\widetilde u\, )+C_5(\,\widetilde v\, )\,
\mod{\bf Z}\).    \hfill\qed

  \section{Abelian extension of \(\Omega^3(SU(n))\)}

 \subsection{ Smooth mappings from \(T=S^3\times [0,1]\) to \(SU(n)\)}
 
Let \(G=SU(n)\), \(n\geq 2\).    Let  \(\Omega^3G\) be the set of smooth mappings from \(S^3\) to \(G=SU(n)\) that
are based at some point \(p_0\in S^3\).        
The 
mapping degree of a \(g\in \Omega^3G\) is given by
\begin{equation}
\deg\,g=\frac{i}{24\pi^2}\int_{S^3}\,
tr(dg\,g^{-1})^3,\end{equation}
where \(dgg^{-1}\) is the pullback of the Maurer-Cartan 1-form on \(G\) by the evaluation map \(ev:\,S^3\times \Omega^3G\longrightarrow G\).   
It satisfies the relation
\begin{equation}\label{deg}
\deg\, ( g_1\cdot g_2)\,=\,\deg g_1+\deg g_2.
\end{equation}
\(g_1\) and \(g_2\) are homotopic if and only if \(\deg g_1=\deg g_2\).    
  $\Omega^3G$ is not connected and is
divided into  denumerable sectors labelled by the mapping degree:
 \begin{eqnarray}
 \Omega^3G&=&\bigcup\,\Omega^3_kG,\\[0.2cm]
 \Omega^3_kG&=&\{g\in \Omega^3G;\quad {\rm deg}\,g=k\,\}.\nonumber
 \end{eqnarray}
  We choose, for each \(k\in {\bf Z}\),   \(g_k\in  \Omega^3_kG\) such that the set  \(\{g_k\}_{k\in {\bf Z}}\) 
is closed under multiplication and that 
\(g_0\equiv 1\) represents  the unit element in \(\pi_3G\).    For example, for \(G=SU(2)\), we may take 
typical instantons
\begin{equation}\label{soliton}
g_k({\bf p})=\left(\begin{array}{cc}s+ir&-q+ip\\[0.2cm]
q+ip&s-ir\end{array}\right)^k,\quad {\bf p}=(p,q,r,s)\in S^3\subset \mathbf{R}^4.
\end{equation}
Here \(g_0=1\) and \(g_{(-k)}\)  is the inverse of \(g_{k}\).    \(g_k( p_0)=1\) at \( p_0=(0,0,0,1)\in S^3\).
From (\ref{deg}) we see that if \(g\in \Omega^3_kG\) then \(g_{(-k)}g\in \Omega^3_0G\).   

Let \(T=S^3\times [0,1]\).    We define the following space of mappings from \(T\) to \(G\).
\begin{eqnarray}\label{K}
K^G&=&\bigcup_{k\in {\bf Z}}\,(K^G)^k\\[0.2cm]
(K^G)^k&=&\left\{u\,:\, S^3\times [0.1] \longrightarrow G\,;\,
\begin{array}{rcl}\frac{\partial u}{\partial t}\vert_{t=0}&=&\frac{\partial u}{\partial t}\vert_{t=1}=0\\[0.2cm]
u(\cdot,0)&=&g_k(\cdot) \,,\\[0.2cm]
u( p_0,\cdot)&=&1 \end{array}\right\}.\nonumber
\end{eqnarray}    
Let 
\begin{equation}\label{J}
J_0^G=\left\{u\in K^G:\, u(\cdot,0)=u(\cdot,1)=1 \right\}\,.
\end{equation}

\begin{lem}\label{u=v}
\begin{enumerate}
\item
Let \(u, v\in K^G\).   If \( u(\cdot,1)=v(\cdot,1)\) then  \( u(\cdot,0)=v(\cdot,0)=g_k(\cdot)\)  for a \(k\in {\bf Z}\).
\item
\(J_0^G\) is the kernel of the map
\begin{equation}
K^G\ni u\,\longrightarrow\, u(\cdot,1)\in \Omega^3G,
\end{equation}
\item
\begin{equation}
\Omega^3G\,=\,K^G/J_0^G.    
\end{equation}
\end{enumerate}
\end{lem}
{\it Proof}

We shall abbreviate \(K^G\), \((K^G)^k\) and \(J^G_0\) to \(K\), \(K^k\) and \(J_0\) respectively.   
If \(u,\,v\in K\) satisfies \(u(\cdot,1)=v(\cdot,1)\), then \(\deg u(\cdot,0)=\deg u(\cdot,1)=\deg v(\cdot,1)=\deg v(\cdot,0)\) which is, say, equal to \(k\).    Then \(u,v\in K^k\) and \( u(\cdot,0)=v(\cdot,0)=g_k(\cdot)\).  
In the same way, if \(u(\cdot ,1)=1\) then \(u(\cdot,0)=1\) and \(u\in J_0\).   This proves assertions 1 and 2.    Now take a \(g\in \Omega^3G\)  such that \(\deg g=k\),   Then  \(g_{(-k)}g\in \Omega^3_0G\) and there is a \(u\in K^0\) such that 
\(u(\cdot,1)=g_{(-k)}g\).   Put  \(v(x,t)=g_k(x)u(x,t)\).   Then \(v\) extends \(g\in \Omega^3_kG\) to \(TG\) and \(v\in K^k\).    This proves the assertion 3.
\hfill\qed

\vspace{0.2cm}

Let  \(T=S^3\times [0,1]\) and \(M=S^3\times S^1\) as before.   We look on \(MG\) as the subset of \(TG\) defined by 
\(\{u\in K^G;\,u(\cdot,0)=u(\cdot,1)\}\).    
Then we have
\begin{equation}\label{JsubMG}
J_0^G=(K^G)^0\cap MG\subset MG\subset K^G.\end{equation}

\subsection{Dual of the space of connections}

Let \({\cal A}_3={\cal A}_3(G)\) be the space of connections on the bundle \(P=S^3\times G\).    Th space  \({\cal A}_3\) is an affine space modeled on the vector space \({\cal E}^1(S^3,\,Lie\,G)\) of \(Lie\,G\)-valued 1-forms on \(S^3\).       
Let \(Map({\cal A}_3,U(1))\) be the group of smooth mappings 
\(\lambda:{\cal A}_3\longrightarrow U(1)\).    The multiplication \(\lambda\cdot\mu\) is defined as the product in \(U(1)\) of their value;  \((\lambda\cdot\mu)(A)=\lambda(A)\cdot\mu(A)\).   
 
 Let \({\cal A}_3^{\ast}\) be the affine dual of \({\cal A}_3\) that is defined by
 \begin{equation}\label{dualspace}
 {\cal A}_3^{\ast}=\left\{\varphi\in  Map({\cal A}_3,\mathbf{ R});\,
 \begin{array}{c} \varphi(pA_1+(1-p)A_2)=p\varphi(A_1)+(1-p)\varphi(A_2),\\[0.2cm]
\, \mbox{ for }\,A_1,A_2\in {\cal A}_3,\quad  p\in\mathbf{ R}\,.\end{array}
\right\}.\\[0.1cm]
 \end{equation}

The directional derivative of a \(\varphi\in Map({\cal A}_3,\mathbf{ R})\) at \(A\) is by definition 
\begin{equation}
(D_A\varphi)a=\lim_{t\longrightarrow 0}\frac{1}{t}(\varphi(A+ta)-\varphi(A)),\quad \forall a\in {\cal E}^1(S^3,\,Lie\,G).\label{der}\end{equation}
If \(\varphi\in {\cal A}_3^{\ast}\), then 
 \begin{equation}
 \varphi(A+a)=\varphi(A)+(D_A\varphi)\, a,\quad \forall a\in{\cal E}^1(S^3,Lie\,G).
 \end{equation}
  \({\cal A}_3^{\ast}\) has a natural vector space structure and if  \(A_0\in {\cal A}_3\) is fixed it is isomorphic to 
  the dual vector space \({\cal E}^1(S^3,Lie\,G)^{\ast}\) by  the correspondence: 
  \[{\cal A}_3^{\ast}\ni\, \varphi \longrightarrow\,D_{A_0}\varphi\in{\cal E}^1(S^3,Lie\,G)^{\ast}.\]
  
   \(\Omega^3G\) acts on \({\cal A}_3^{\ast}\) by
 \begin{equation}
 (f\cdot \varphi)(A)=\varphi(f\cdot A), \quad\mbox{ for } \varphi\in {\cal A}_3^{\ast}\,,\,\quad f\in \Omega^3G,
\end{equation}
 where \[
f\cdot A=f^{-1}Af+f^{-1}df.\]
We have
\begin{equation}
(D_A(f\cdot\varphi ))a=(D_{f\cdot A}\,\varphi) \,Ad_f a,\quad\,\forall a\in {\cal E}^1(S^3,\,Lie\,G).
\end{equation}

  \begin{defn}\label{character}
  \begin{equation}
\exp 2\pi i\,{\cal A}_3^{\ast}=\left\{ \exp 2\pi i\varphi\in Map({\cal A}_3,U(1))\,;\quad \varphi\in {\cal A}_3^{\ast} \right\}
\end{equation}
 \end{defn}
 
 \(\exp 2\pi i\,{\cal A}_3^{\ast}\) is an abelian subgroup of \(Map({\cal A}_3,U(1))\).

 \subsection{Abelian extension of \(\Omega^3(SU(n))\), \(n\geq 3\)}
 
 In this section we study the abelian extension of \(\Omega^3 G\) for the case \(G=SU(n)\), \(n\geq 3\). .     In ~\cite{Mi, Mic}  Mickelsson described the extension of \(\Omega^3_0 G=(K^G)^0/J^G_0\,\) by the abelian group \(Map({\cal A}_3,U(1))\).    We shall extend his argument to the case 
 \(\Omega^3G=(K^G)/J^G_0\) and the extension is given by the abelian group  \(\exp 2\pi i\,{\cal A}_3^{\ast}\).     For \(G=SU(n)\) with  \(n\geq 3\) we have \(\pi_4(G)=0\), so  \(C_5(g)\) for \(g\in J_0\subset MG\) is well defined, (\ref{JsubMG}) and (\ref{C5}).      
  In particular 
Polyakov-Wiegmann formula of Lemma \ref{Poly} is valid  for \(f,\,g\in
J_0\).
\begin{equation}
C_5(fg)=C_5(f)+C_5(g)+ \beta_T(f,g)\quad
\mod{\bf Z}.\label{PolyJ0}\end{equation}
\begin{defn}
For \(f\in K\) and \(g\in J_0\), we put
\begin{equation}
\alpha_T(f,g)=\beta_T(f,g)+C_5(g).\label{alpha}
\end{equation}
\end{defn}
\begin{lem}\label{cocycle} 
Let \(G=SU(n)\), \(n\geq 3\), and \(T=S^3\times [0,1]\).     Let \(f, g,h\in
TG\) be such that \(f(\cdot,1)=g(\cdot,1)=h(\cdot,1)\).    Then we have
\begin{equation}
\alpha_T(f,f^{-1}g)+\alpha_T(g,g^{-1}h)=\alpha_T(f,f^{-1}h) \mod {\bf Z}.
\end{equation}
\end{lem}
In fact,  since \(f^{-1}g\,,\,g^{-1}h,\,f^{-1}h\in J_0\) we have from (\ref{PolyJ0}) 
\[C_5(f^{-1}g)+C_5(g^{-1}h)+\beta_T(f^{-1}g,\,g^{-1}h)=C_5(f^{-1}h).\]
On the other hand (\ref{beta}) implies
\[\beta_T(f,f^{-1}g)+\beta_T(g,g^{-1}h)=\beta_T(f^{-1}g,g^{-1}h)+\beta_T(f,f^{-1}h) .\]
   Hence Definition (\ref{alpha}) of \(\alpha_T\) yields the desired equation.   \hfill\qed

\begin{defn}
For  \(f, \,g\in K\) such that \(f(\cdot,1) =g(\cdot,1)\), we put 
\begin{equation}
\chi_T(f,g)=\exp 2\pi i\, \alpha_T(f, f^{-1}g).\label{chi}
\end{equation}
\end{defn}
From Lemma \ref{cocycle} we have 
\begin{equation}
\chi_T(f,g)\chi_T(g,h)=\chi_T(f,h),\end{equation}
for \( f,g,h\in K\) such that  \(f(\cdot,1) =g(\cdot,1)=h(\cdot,1)\).

  We define the right action of \(J_0\) on the product set \(K \times Map({\cal A}_3,U(1))\) by 
\begin{equation}
g\cdot\,(f\,,\lambda\,)=\,\left(\,f\,g\,,\,\lambda(\cdot) \chi_T(f, fg)\,\right),\qquad g\in J_0.
 \end{equation}
 Note that \(J_0\) acts trivially on \({\cal A}_3\), and on \(Map({\cal A}_3,U(1))\).   
We consider the quotient space:
\begin{equation}
\widehat{\Omega G}=K\times \exp 2\pi i{\cal A}_3^{\ast}/J_0\,.\label{extension}
\end{equation}
It is the quotient by the equivalence relation
\[(f,\lambda)\sim (g,\mu)\,\Longleftrightarrow\,\left\{\begin{array}{lcr}
 g(\cdot,1) &=& f(\cdot,1)\qquad\\[0.2cm]
\mu(\cdot)&=&\lambda(\cdot)\,\chi_T(f,g)\,.
\end{array}\right.  \]
The equivalence class of \((f,\lambda)\) is denoted by \([f,\lambda]\), and  
the projection \(\widehat \pi: \widehat{\Omega G}\longrightarrow \Omega^3G\) is given by \(\widehat \pi([f\,,\,\lambda])=f(\cdot,1)\).   Then \(\widehat{\Omega G}\) becomes a principal bundle over \(\Omega^3G\) with the structure group \(\exp 2\pi i{\cal A}_3^{\ast}\,\).  Here the \(U(1)\) valued transition function \(\chi_T(f,g)\) is considered as a constant map in \(\exp 2\pi i{\cal A}_3^{\ast}\).    

\begin{thm}
Let \(G=SU(n)\) with \(n\geq 3\).   Then 
\(\widehat{\Omega G}\) gives a Lie group extension of \(\Omega^3G\) by the abelian group 
 \(\exp 2\pi i{\cal A}_3^{\ast}\,\).
\end{thm}

{\it Proof}

We shall endow \(\widehat{\Omega G}\)  with a group structure by  Mickelsson's 
2-cocycle (\ref{Mic3}):
\begin{eqnarray}
 \gamma_{T}(f,g\,; A)&=&
-\frac{i}{24\pi^3}\int_{T}(\delta c^{1,1})(\,f,\,g\,;\,A) \nonumber
\\[0.2cm]
&=& \frac{i}{24\pi^3}\int_{\partial T}c^{2,0}(f,g\,;\,A)+
\beta_{T}(f,g)\,. \label{Mic3n}
\end{eqnarray}
It holds that \(\gamma_T(f,g;\cdot)\in {\cal A}_3^{\ast}\) and \(\exp 2\pi i\gamma_{T}(f,g\,;\cdot) \in \exp 2\pi i{\cal A}_3^{\ast}\).   
We define the product on \(K\times \exp 2\pi i{\cal A}_3^{\ast}\) by
\begin{equation}\label{gammulti}
(f,\lambda)\ast(g,\mu)=\left(\, fg
,\,\lambda(\cdot)\mu_{f}(\cdot)\exp 2\pi i\gamma_{T}(f,g\,;\,\cdot)\,\,\right).
\end{equation}
 The associativity of the product follows from the cocycle condition (\ref{gamma}) of \(\gamma_T\).    
We shall verify that the multiplication rule on \(K\times \exp 2\pi i{\cal A}_3^{\ast}\) descends to that on \(\widehat{\Omega G}\).   
Let  \((f,\lambda)\), \((g,\mu)\) and \((h,\nu)\in K\times \exp 2\pi i{\cal A}_3^{\ast}\) and suppose  \((f,\lambda) \sim (g,\mu)\).    Since \(g(\cdot,1) =f(\cdot,1) \)  it holds that \((fh)(\cdot,1) =(gh)(\cdot,1)\), so Lemma \ref{u=v} implies \( fh=gh\) on \(\partial T=S^3\times\{0,1\}\), and we have 
\begin{equation}
\int_{\partial T}c^{2,0}(g,h\,;\,A)=\int_{\partial T}c^{2,0}(f,h\,;\,A).
\end{equation}
From (\ref{beta}),~(\ref{beta2})
 and Lemma~\ref{MiC5} we have 
 \begin{equation}\label{multi}
\alpha_{ T}(f,f^{-1}g)+\beta_{ T}(g,h)=  \alpha_{ T}(fh,h^{-1}f^{-1}gh)+ \beta_{T}(f,h)\,.
\end{equation}
These two equations imply 
\begin{equation*}
\exp 2\pi i\gamma_T(g,h:\cdot)\chi_T(f,g)=\exp 2\pi i\gamma_T(f,h:\cdot)\chi_T(fh,gh). 
\end{equation*}
Now from the assumption we have \(\mu(\cdot)=\lambda(\cdot)\chi_T(f,g)\) and \(\nu_g(\cdot)=\nu_f(\cdot)\).    Hence
\begin{equation*}
\mu(\cdot)\nu_g(\cdot)\exp 2\pi i\gamma_T(g,h:\cdot)=\lambda(\cdot)\nu_f(\cdot)\exp 2\pi i\gamma_T(f,h:\cdot)\chi_T(fh,gh).
\end{equation*}
Therefore  \((f,\lambda)\ast (h,\nu)\sim(g,\mu)\ast (h,\nu)\).      
Next we suppose  \((g,\mu) \sim (h,\nu)\).    Then \(\nu_f(\cdot)=\mu_f(\cdot)\chi_T(g,h)\).    By the same calculations as above we have
\begin{equation*}
\exp 2\pi i\gamma_T(f,h:\cdot)\chi_T(g,h)=\exp 2\pi i\gamma_T(f,g:\cdot)\chi_T(fg,fh).
\end{equation*}
   Hence
\begin{equation}
\lambda(\cdot)\nu_f(\cdot)\exp 2\pi i\gamma_T(f,h:\cdot)=\lambda(\cdot)\mu_f(\cdot)\exp 2\pi i\gamma_T(f,g:\cdot)\chi_T(fg,fh).\label{equiv}
\end{equation}
This implies \((f,\lambda)\ast (g,\mu)\sim (f,\lambda)\ast (h,\nu)\).      Therefore 
\(\widehat{\Omega G}\) inherits the group structure.    The unit of \(\widehat{\Omega G}\) is 
the equivalence class \([1,1]\) consisting of \(1\in K^0\) and the constant map \(1\in \exp 2\pi i{\cal A}_3^{\ast}\).    The inverse of  \([f,\lambda]\) is \([f^{-1},(\lambda_{f^{-1}})^{-1}]\), 
where \((\lambda_{f^{-1}})^{-1}(A)=(\lambda(f^{-1}\cdot A))^{-1}\), and we used the fact \(\gamma_T(f,f^{-1};\cdot)=0\) of Lemma~\ref{beta2}.
 The group \(\,\exp 2\pi i{\cal A}^{\ast}_3\,\) is embedded as a normal subgroup of \(\widehat{\Omega G}\) by the map \(\lambda\longrightarrow\,[1,\lambda] \in \widehat{\Omega G}\) .   Thus  \(\widehat{\Omega G}\)  is  an extension of \(\Omega^3 G\) by the abelian group  \(\,\exp 2\pi i{\cal A}_3^{\ast}\,\).       \hfill\qed

\subsection{Extension of the embedded  \(\Omega^3(SU(2))\)}

We shall denote \(G=SU(2)\) and \(G^{\prime}=SU(3)\).     \(\Omega^3G\) being embedded in \(\Omega^3G^{\prime}\) we may think that the restriction to \(\Omega^3G\) of the group extension \(\widehat{\Omega G^{\prime}}\)  of  \(\Omega^3G^{\prime}\) 
yields a group extension of \(\Omega^3 G\).    Let  \(\widehat\pi:\,\widehat{\Omega G^{\prime}} \longrightarrow \Omega^3 G^{\prime}\)  be the extension of \(\Omega^3 G^{\prime}\) discussed in  subsection 2.3.    
 Let  \(K\) and \(J_0\) ( respectively  \(K^{\prime}\) and \(J_0^{\prime}\) ) be the spaces defined by the formulas  (\ref{K}) and (\ref{J})  for the group \(G=SU(2)\) ( respectively for the group \(G^{\prime}=SU(3)\)).     
      Let 
   \begin{eqnarray*}
  \widetilde K&=&\left\{\widetilde u=\left(\begin{array}{cc}u&0\\[0.2cm]0&1\end{array}\right):\, S^3\times [0.1] \longrightarrow G^{\prime}\,;\, u\in K\,\right\}\subset K^{\prime}.\\[0.2cm]
  \widetilde J_0&=&\left\{\widetilde u=\left(\begin{array}{cc}u&0\\[0.2cm]0&1\end{array}\right)\in \widetilde K\,;\,u\in J_0\right\}\subset J_0^{\prime}.
      \end{eqnarray*}
      We have
  \begin{equation*}
\widehat{\Omega G^{\prime}}=K^{\prime}\times \exp 2\pi i{\cal A}_3^{\ast}/J^{\prime}_0\,,\qquad (\ref{extension})
\end{equation*}
with the transition function given by (\ref{chi});
\[\chi_T(f^{\prime},g^{\prime})=\exp 2\pi i\alpha_T(f^{\prime},(f^{\prime})^{-1}g^{\prime}),\quad f^{\prime},g^{\prime}\in K^{\prime}.\] 
The restriction to \(  \widetilde K\)  becomes
 \begin{equation}
 \chi_T(\widetilde f, \widetilde g)=
 \exp 2\pi i\,C_5(\widetilde f^{-1}\widetilde g)=\epsilon (f^{-1}g)\,,\quad  f, g\in K.   
 \end{equation}
It defines a  \(U(1)\)-principal  bundle over \(\Omega^3G\).   In fact, 
 we define the right action of \(\widetilde J_0\) on the space \(\widetilde K\times U(1)\) by 
  \begin{equation}
  \widetilde g\cdot (\widetilde f,c)= (\,\widetilde f\,\widetilde g\,, \,\epsilon(\widetilde g)\,c\,) \quad\mbox{ for  } (\widetilde f,c )\in \widetilde K\times U(1) \,\mbox{ and } \widetilde g\in \widetilde J_0.\label{Jaction}
  \end{equation}
Then we have the \(U(1)\)-principal  bundle 
\begin{equation}
\widehat{\Omega G}=\widetilde K\times U(1)/\widetilde J_0\,\longrightarrow\,\Omega^3G.
\end{equation}
\(\widehat{\Omega G}\) is reduced to a \({\bf Z}_2\)-principal bundle.   From the discussion preceding 
  Definition~{\ref{epsilon} we see that  \(\widehat{\Omega G}\) represents the action of the homotopy group \(\pi_4(SU(2))\) on \(\Omega^3G\).   We note that the Mickelsson's 2-cocycle (\ref{Mic3}) vanishes from Proposition \ref{nul} so that the group extension is algebraically trivial though it is not topologically trivial.

 \begin{prop}
 Let \(G=SU(2)\).   
 \begin{enumerate}
 \item
 There is a Lie group extension \(\widehat{\Omega G}\) of \(\Omega^3G\) by \({\bf Z}_2\).   
 \item
The restriction of the bundle \(\widehat{\Omega G}\) to \(\Omega^3_{ev}G=\,\cup_{k:\,even}\Omega^3_kG\) is isomorphic to  \(\Omega^3_{ev}G\), and the restriction to  \(\Omega_{0}G\)  is isomorphic to  \(\Omega^3_{0}G\).
\end{enumerate}
\end{prop} 
In fact, let 
 \(f,\,g\in K\) be such that \(f(\cdot,1)=g(\cdot,1)\).   If \(f(\cdot,1)=g(\cdot,1)\) is in \(\Omega^3_{ev}G\)  then  \(f^{-1}g\) is in the trivial homotopy class of \(\pi_4G\) and \(\epsilon(f^{-1}g)=1\).      \hfill\qed

\section{Adjoint representations of \(\widehat{\Omega_0G}\) }

\subsection{Lie algebra extensions of \(S^3(Lie G)\) } 

Let \(G=SU(n)\) with \(n\geq 3\).     Let \(\widehat{\Omega_0G}\) be the Lie group extension of \(\Omega^3_0G\) :
 \[\widehat{\Omega_0G}=
K^0\times\exp 2\pi i{\cal A}_3^{\ast}/J_0.   \]
We shall study the Lie algebra of the Lie group extension \(\widehat{\Omega_0G}\).
   In the following elements of the dual space \({\cal A}_3^{\ast}\) are denoted by \(l,\,m\,\), and an element of \(\widehat{\Omega_0G}\) is denoted by the pair  \((f,\lambda)\) with \(f\in K^0\) and \(\lambda(\cdot)=\,e^{2\pi i \,l(\cdot)}\,\in \exp 2\pi i{\cal A}_3^{\ast}\).      
We note that the 2-cocycle (\ref{Mic3n}) becomes 
\begin{equation}
\gamma_T(f,g;A)= \frac{i}{24\pi^3}\int_{S^3}c^{2,0}(f,g\,;\,A)+\beta_{T}(f,g),
\end{equation}
because the integral on the part of the boundary  \(S^3\times \{0\}\subset \partial T\) vanishes by virtue of the facts  \(\Omega_0G=K^0/J_0\) and \(K^0\simeq D^4G\) .   
 \(S^3(Lie\,G)\) denotes the Lie algebra of the based mappings from \(S^3\) to \(Lie \,G\).

 \begin{thm}~~
The Lie algebra of  \(\widehat{\Omega_0G}\) is 
 given by the 
vector space \(S^3(Lie \,G)\oplus {\cal A}_3^{\ast}\) endowed with the 
commutation relation 
\begin{equation}\label{Liealgext}
\left[\,(\xi,\,l),\,(\eta,\,m)\,\right]=\left(\,[\xi,\eta]\,,\,D_Am(d_A\xi)\,-\,D_Al(d_A\eta)\,+\,i\omega(\xi,\eta\,;\,A\,)\,\right),
\end{equation}
where
\(D_Am\) is the derivation of \(m\),~(\ref{der}), 
and
\begin{equation}
\omega(\xi,\eta\,;A)=-\frac{1}{24\pi^3}\int_{S^3}\,tr(d\xi d\eta-d\eta d\xi)A\,.\label{Liebra}
\end{equation}
The Lie algebra \(\left(\,S^3(Lie \,G)\oplus{\cal A}_3^{\ast},\,[\,\cdot,\cdot\,]\,\right)\) becomes the Lie algebra extension of \(S^3(Lie\,G)\).   
\end{thm}

{\it Proof}

Recall that the group multiplication law for \(\widehat{\Omega SU(n)}\), \(n\geq 3\) is defined in  (\ref{gammulti}).    We have 
\begin{eqnarray*}
&&(e^{s\xi},1)(e^{t\eta},1)(e^{-s\xi},1)(e^{-t\eta},1)=\left(e^{s\xi}e^{t\eta}e^{-s\xi}e^{-t\eta}, \,\exp 2 \psi(s,t;\cdot)\right),
\\[0.2cm]
&&\psi(s,t;A)=
\,\gamma_T(e^{s\xi},e^{t\eta}\,;A)+\gamma_T(e^{s\xi}e^{t\eta},e^{-s\xi}\,;A)+\gamma_T(e^{s\xi}e^{t\eta}e^{-s\xi},e^{-t\eta}\,;A),
\end{eqnarray*}
where \(1\) means the constant map \({\cal A}_3\ni A\longrightarrow 1\in U(1)\).
Since 
\begin{eqnarray*}
&&\frac{d}{ds}\vert_{s=0}\frac{d}{dt}\vert_{t=0}\,c^{2,1}(e^{s\xi},e^{t\eta})=0,\\[0.2cm]
&&
\frac{d}{ds}\vert_{s=0}\frac{d}{dt}\vert_{t=0}\,c^{2,0}(e^{s\xi},e^{t\eta};A)=
\frac{d}{ds}\vert_{s=0}\frac{d}{dt}\vert_{t=0}\,c^{2,0}(e^{s\xi}e^{t\eta},e^{-s\xi};A)
= e^{2,0}(\xi,\eta;A),\\[0.2cm]
&&\frac{d}{ds}\vert_{s=0}\frac{d}{dt}\vert_{t=0}\,c^{2,0}(e^{s\xi}e^{t\eta}e^{-s\xi},e^{-t\eta};A)=0,
\end{eqnarray*}
we have
\begin{eqnarray*}
2  \frac{d}{ds}\vert_{s=0}\frac{d}{dt}\vert_{t=0}\,\psi(s,t;A)&=& -\frac{i}{12\pi^3}\int_{S^3}\,e^{2.0}(\xi,\eta\,;\,A)\\[0.2cm]
&=&-\frac{i}{24\pi^3}\int_{S^3}\,tr[(d\xi d\eta-d\eta d\xi)A].
\end{eqnarray*}
Therefore 
\[\frac{d}{ds}\vert_{s=0}\frac{d}{dt}\vert_{t=0}\,(e^{s\xi},1)(e^{t\eta},1)(e^{-s\xi},1)(e^{-t\eta},1)=\left([\xi,\eta],\,-i\omega(\xi,\eta\,;\,\cdot)\right).\]
The commutation relation (\ref{Liealgext}) follows from this.

\subsection{ Adjoint orbits of \(\widehat{\Omega_0 G}\)}

Let \(\widehat{S^3(Lie\,G)}=S^3(Lie\,G)\oplus {\cal A}^{\ast}_3\) be the Lie algebra extension of \(S^3(Lie\,G)\).   

\begin{prop}
The adjoint action of a \((g,\nu)\in \widehat{\Omega_0G}\) on \(Lie(\,\widehat{\Omega_0G}\,)\,=\widehat{S^3( Lie\,G)}\)  is described as follows.
\begin{equation*}
Ad_{(g,\nu)}\,(\xi,l\,)\vert_A=\frac{d}{ds}\vert_{s=0}\,Ad_{(g,\nu)}\,(\,e^{s\xi},\,e^{\,s\,l(\cdot)})\vert_A
=\left(\,Ad_g\xi\,, \,O(\xi,l;g,\nu)\right),\end{equation*}
where \(\nu=e^{m(\cdot)}\) with  \(m\in {\cal A}_3^{\ast}\)  and 
\begin{eqnarray*}
\,O(\xi,l;g,e^{m(A)})\,&=&g\cdot l(A)- \, D_Am(\,g\,d_{g\,\cdot A}\xi\,g^{-1})\, \\[0.2cm]
&&+\frac{1}{12\pi^2}\int_{S^3}\,tr\left(\,(g^{-1}dg\,d\xi-d\xi \,g^{-1}dg)g^{-1}Ag\, -d\xi(g^{-1}dg)^3 \, \,\right)\\[0.2cm]
&&+\frac{1}{24\pi^2}\int_{ S^3}\,tr\,\left(\,[\,g^{-1}dg,\,[g^{-1}dg,\,\xi\,]\,]\,g^{-1}Ag\,\right)\,,
\end{eqnarray*}
\end{prop}
for \((\xi,l)\in\,\widehat{S^3(Lie G)}\).

{\it Proof}

The adjoint action of \( \widehat{\Omega_0G}\) is given by 
\begin{equation*}
Ad_{(g,\nu)}\,(f,\lambda)=\left(\,gfg^{-1}\,,\,\nu(\cdot)\lambda_g(\cdot)(\nu_{gfg^{-1}}(\cdot))^{-1}
\exp 2\pi i\{\gamma_T(g,f;\cdot)+\gamma_T(gf,g^{-1};\,\cdot\,)\}\,\right),
\end{equation*}
for \((f,\lambda),\,(g,\nu )\in \widehat{\Omega_0G}\).   
If   \((f,\lambda)=\left(e^{s\xi} ,\,e^{\,s\,l(\cdot)}\right)\) with  \(\,\xi\in S^3(Lie\,G)\),  \(l\in {\cal A}_3^{\ast}\), and if 
\((g, \nu)=\left(g, e^{\,m(\cdot)}\right) \) with \(m \in {\cal A}_3^{\ast}\), we have 
\begin{eqnarray*}
\frac{d}{ds}\vert_{s=0}\,e^{sl(g\cdot A)}&=&l(g\cdot A),\\[0.2cm]
\frac{d}{ds}\vert_{s=0}\,\left(\nu( (ge^{s\xi}g^{-1})\cdot A)\right)^{-1}
&=& -(D_Am)(\,g\,d_{g\cdot A}\xi\,g^{-1})
,\\[0.2cm]
\frac{d}{ds}\vert_{s=0}\,\gamma_T(g,e^{s\xi};A)&=&
\frac{i}{48\pi^3}\int_{D^4}\,tr[\,d\xi(g^{-1}dg)^3]\\[0.2cm]
\quad&& -\frac{i}{48\pi^3}\int_{S^3}\,tr[\,(g^{-1}dg\,d\xi-d\xi \,g^{-1}dg)g^{-1}Ag\,]
,\\[0.2cm]
\frac{d}{ds}\vert_{s=0}\,\gamma_T(ge^{s\xi},g^{-1};A)&
=&\frac{i}{48\pi^3}\int_{D^4}\,tr[\,(d_{g^{-1}dg}\xi)(g^{-1}dg)^3]\\[0.2cm]
&&-\frac{i}{48\pi^3}\int_{S^3}\,tr[\,\left((g^{-1}dg)\,d_{g^{-1}dg}\xi-d_{g^{-1}dg}\xi\,( g^{-1}dg)\right)g^{-1}Ag\,].
\end{eqnarray*}
Hence the adjoint representation of \((g,\nu)\) on the Lie algebra \(\widehat{S^3( Lie\,G)}\) is given by 
\begin{equation*}
Ad_{(g,\nu)}\,(\xi,l\,)\vert_A=\frac{d}{ds}\vert_{s=0}\,Ad_{(g,\nu)}\,(\,e^{s\xi},\,e^{\,s\,l(\cdot)})\vert_A
=(Ad_g\xi,O(\xi,l;g,\nu)\,),
\end{equation*}
with
\begin{eqnarray*}
O(\xi,l;g,\nu)\vert_A&=&
g\cdot l(A)- \, D_Am(\,g\,d_{g\cdot A}\xi\,g^{-1})\\[0.2cm]
&&+\frac{1}{12\pi^2}\int_{S^3}\,tr\left(\,(g^{-1}dg\,d\xi-d\xi \,g^{-1}dg)g^{-1}Ag\, -d\xi(g^{-1}dg)^3 \,  \right)\\[0.2cm]
&&+\frac{1}{24\pi^2}\int_{ S^3}\,tr\,\left(\,[\,g^{-1}dg,\,[g^{-1}dg,\,\xi]\,]\,g^{-1}Ag\,\right).
\end{eqnarray*}

 \hfill\qed

\begin{prop}    
The adjoint action \( ad\in End(\, \widehat{S^3( Lie\,G)}\,)\) becomes 
\begin{equation}
 ad_{ (\eta, m) }\, (\xi,l\,) \vert_A=\left(\,[\xi,\eta]\,,\,
D_Al(d_A\eta)-\,D_Am(d_A\xi)-\,i\,\omega(\xi,\eta\,;A)\,\right),
\end{equation}
where 
\[\omega(\xi,\eta\,;A)=-\frac{1}{24\pi^3}
\int_{S^3}\,tr [(d\xi\,d\eta -d\eta\,d\xi)A\,].\]
It coincides with the Lie bracket of \(\widehat{S^3 (Lie\,G})\), (\ref{Liebra}).
\end{prop}

In fact, if we let \((g,\nu)=(e^{t\eta},e^{tm})\) in \(Ad_{(g,\nu)}(e^{s\xi},e^{sl})\), 
\begin{eqnarray*}
\frac{d}{ds}\vert_{s=0}\frac{d}{dt}\vert_{t=0}\,e^{sl(e^{t\eta}\cdot A)}&=&\frac{d}{dt}\vert_{t=0}l(e^{t\eta}\cdot A)=D_Al(d_A\eta),\\[0.2cm]
\frac{d}{ds}\vert_{s=0}\frac{d}{dt}\vert_{t=0}\left(\nu( (ge^{s\xi}g^{-1})\cdot A)\right)^{-1}
&=& -\frac{d}{dt}\vert_{t=0}(D_A(tm))(e^{t\eta}\,d_{e^{t\eta}\cdot A}\xi\, e^{-t\eta})\\[0.2cm]
&=&-D_Am(d_A\xi).
\end{eqnarray*}
 \hfill\qed

  \end{document}